\documentclass{article}
\usepackage{graphicx} 
\usepackage{xurl} 
\usepackage{hyperref}
\usepackage{xcolor}
\usepackage{caption}
\usepackage{float}
\usepackage{subcaption}
\usepackage{pifont}
\usepackage{booktabs}
\usepackage{svg}
\PassOptionsToPackage{FIGTOPCAP}{subfigure} 
\usepackage{siunitx}
\usepackage{amsmath}
\usepackage{amssymb}
\usepackage{textcomp,gensymb}	
\usepackage{amsmath}

\usepackage{comment}
\newcommand{\COtwo}{CO\textsubscript{2} }
\begin{document}


\begin{titlepage}






\hspace{-.55cm}{\fontfamily{pbk}\selectfont
\large From carbon management strategies to implementation: \\}
{\fontfamily{pbk}\selectfont \noindent Modeling and physical simulation of CO$_2$ pipeline infrastructure - \\a case study for Germany}

\vspace*{1cm} 		
{\small\noindent Mehrnaz Anvari$^{1}$, Marius Neuwirth$^{2}$, Okan Akca$^{1,3}$, Luna Lütz$^{2,3}$, Simon Lukas Bussmann$^{2,4}$, Tobias Fleiter$^{2}$, Bernhard Klaassen$^{5,*}$}

\vspace*{1cm}

\hspace{-.55cm}\textbf{Abstract}\\
Carbon capture and storage or utilization (CCUS) will play an important role to achieve climate neutrality in many economies. 
Pipelines are widely regarded as the most efficient means of CO$_2$ transport; however, they are currently non-existent. Policy-makers and companies need to develop large-scale infrastructure  under substantial uncertainty. Methods and analyses are needed to support pipeline planning and strategy development. 

This paper presents an integrated method for designing CO$_2$ pipeline networks by combining energy system scenarios with physical network simulation. Using Germany as a case study in a projection to the year 2045, we derive spatially highly resolved CO$_2$ balances to develop a dense-phase CO$_2$ pipeline topology that follows existing gas pipeline corridors. The analyzed system includes existing sites for cement and lime production, waste incineration, carbon users, four coastal CO$_2$ hubs, and border crossing points. We then apply the multiphysical network simulator MYNTS to assess the technical feasibility of this network. We determine pipeline diameters, pump locations, and operating conditions that ensure stable dense-phase transport. The method explicitly accounts for elevation and possible impurities.

The results indicate that a system of about 7000 km pipeline length and a mixed normed diameter of DN700 on main corridors and of DN500/DN400 on branches presents a feasible solution to connect most sites. Investment costs for the optimized pipeline system are calculated to be about 17 billion Euros. 
The method provides a reproducible framework and is transferable to other countries and to European scope.

\vspace*{0.5cm}
\noindent Keywords: Carbon capture and storage, modeling, simulation, pipeline transport, dense phase  

\vfill

\hspace{-.55cm}
\begin{tabular}{@{}lp{\linewidth}@{}}
$^1$ & Fraunhofer Institute for Algorithms and Scientific Computing SCAI \\
$^2$ & Fraunhofer Institute for Systems- \& Innovation Research ISI\\
$^3$ & Technical University of Berlin  \\
$^4$ & Kassel University \\
$^5$ & Fraunhofer Research Institution for Energy Infrastructure and Geotechn. IEG\\
$^*$ & Corresponding author: Bernhard Klaassen, Fraunhofer IZB, Schloss Birlinghoven 1,  53757 Sankt Augsutin, Germany, Tel.: +49 2241 144070, 
e-mail: bernhard.klaassen@ieg.fraunhofer.de
\end{tabular}

\date{January 2026}


\end{titlepage}



For large CO$_2$ volumes the most effective option, pipelines are usually designed to transport CO$_2$ in the dense phase over long distances.
\begin{itemize}
    \item \textbf{Trains:} Suitable for smaller quantities and a flexible solution for individual sites \cite{Becattini.2024}. Rail transport can be deployed immediately and can therefore serve as a bridge until pipeline connections are established.
    \item \textbf{Ships:} For international transport, CO$_2$ can be liquefied and shipped by specialized vessels to locations where it is stored or utilized. Inland waterway vessels can also be used, but they have significantly lower capacity than seagoing ships and depend on navigable water levels \cite{Becattini.2024}.
    \item \textbf{Trucks:} A feasible option only for very small quantities and generally not competitive for large-scale CO$_2$ transport.
\end{itemize}

In this study, we focus on CO$_2$ transport via pipelines. The development of a CO$_2$ system and transport infrastructure is still at the very early stage and most countries do not operate any CO$_2$ pipelines. Several countries, however, have published strategies for the development of CCUS and CO$_2$ infrastructure, including Germany \cite{FederalGovermentGermany.2025, Bundestag.}, France \cite{Ministeredelatransitinecologiqueetsolidaire.2020}, and the Netherlands \cite{MinisterievanEconomischeZaken.2025}. In Germany, the recently adopted Carbon Storage and Transport Act (KSpTG) \cite{FederalGovermentGermany.2025} provides, for the first time, a comprehensive legal framework for both the permanent storage of CO$_2$ (primarily in geological formations in the North Sea) and the permitting and operation of CO$_2$ pipelines. The act classifies CO$_2$ pipelines and storage as being in the overriding public interest and allows the transport and storage of CO$_2$ from all sources except coal-fired power generation. At the same time, the law does not foresee a centrally planned national backbone; instead, the development of CO$_2$ pipeline infrastructure is initially left to private and public project developers, with access and connection rules defined at the network level. In parallel, the federal government is preparing a national carbon management strategy, which is expected to further specify the role of CO$_2$ capture, transport, and storage in achieving Germany’s climate targets \cite{BundesministeriumfurWirtschaftundKlimaschutz.26.02.2024}. Despite these advances, the actual rollout of CO$_2$ pipeline networks has not yet begun. Also at European level,  a clear long-term vision for a coherent European CO$_2$ pipeline system is still lacking.  
In order to facilitate progress with infrastructure planning and strategies, analytical methods need to integrate systems analysis with concrete pipeline planning. In the past, most energy system studies have been very generic with regard to specific CO$_2$ pipeline infrastructure while on the other side, pipeline planning was lacking the system context and the long-term vision.

We aim to fill this gap by developing a method to calculate CO$_2$ pipeline networks for entire regions or countries by combining physical pipeline simulation with energy system modeling, and we apply it to a potential CO$_2$ network for a future climate-neutral energy system in Germany. The method is intended to provide a basis for decision makers in industry and policy to develop a long-term vision and concrete plans for CO$_2$ infrastructure in the future energy system.

There is an increasing body of literature on CO$_2$ pipeline networks in Europe \cite{Morbee.2012, dAmore.2021, Oeuvray.2024} and Germany \cite{Hofmann.2025, Nguyen.2024, Mohr.2024, Yeates.2024, DanielBenrathStefanFlammeSabrinaGlanzFanziskaM.Hoffart.31.08.2020}. These studies differ in their approaches to topology development and spatial resolution. Some assume straight-line connections between sources and sinks \cite{Morbee.2012, dAmore.2021, Nguyen.2024, Hofmann.2025}, others follow existing natural gas transmission pipelines \cite{Mohr.2024, DanielBenrathStefanFlammeSabrinaGlanzFanziskaM.Hoffart.31.08.2020}, and others use more flexible routing methods \cite{Yeates.2024}. The sets of CO$_2$ sources also differ: most studies include sectors similar to those considered here, but often add fossil sources that are excluded in this paper because alternative mitigation options are assumed to be more appropriate for those emissions. The inclusion of basic chemicals via CCU has so far not been explicitly considered in comparable German studies. 

Among the existing work, the VDZ study \cite{VDZ.2024} is noteworthy, as it develops a pipeline topology in Germany for waste incineration, cement, and lime and estimates associated costs, although the underlying calculations are not fully transparent. Another important analysis was developed in the ELEGANCY project by Benrath et al. \cite{DanielBenrathStefanFlammeSabrinaGlanzFanziskaM.Hoffart.31.08.2020}, where a CO$_2$ pipeline network along the gas grid was combined with simplified pressure and pump placement calculations to derive system and transport costs. At the European level, the JRC has conducted a study \cite{Tumara.2024} that uses straight-line connections, aggregates various locations, and includes a ramp-up analysis of CO$_2$ infrastructure. A distinctive feature of our analysis is the explicit coupling of a scenario-based CO$_2$ topology with a physically pipeline simulation. While similar work has recently been published for Austria, developed jointly with a network operator by Schützenhofer et al. \cite{Schützenhofer.10.2024}, including pressure calculations and pump placement, this type of analysis has not yet been carried out for Germany. 

To demonstrate and validate the method, we develop a case study for a climate-neutral energy system in Germany. We first construct a CO$_2$ network topology for Germany, including locations of CO$_2$ sources and sinks and potential transport routes for CO$_2$ from neighboring countries. Based on prevailing uncertainties, we define several scenarios with different assumptions about transit, capture, and utilization volumes. We then apply the scenario with the highest infrastructure requirements to the topology and, in a second step, simulate CO$_2$ flows within the network using a fluid dynamics modeling approach. This enables the optimization of key network parameters, including pipeline diameters, pump placement, and resulting infrastructure costs. The method is generic and can be transferred to other countries or regions as well as adapted scenarios.

The following sections answer these research questions:
\begin{itemize}
    \item What volumes of CO$_2$ will be generated in the sectors under consideration for CO$_2$ transport in Germany in 2045?
    \item How can a dense-phase pipeline backbone - based on these volumes - be designed along the natural gas network?
    \item How do elevation differences, CO$_2$ impurities, temperatures, and pipe diameters affect CO$_2$ transport and the associated costs?
\end{itemize}

\newpage

\section{Method and data}

\subsection{Overview of the method}
Our method combines the broad perspective and long-term vision of energy system modeling with the detailed technical infrastructure planning. From energy system scenarios for a future climate neutral industry we derive specific CO$_2$ infrastructure needs and use the pipeline planning tool MYNTS to calculate system parameters like the diameter of pipelines and the need for pumps and compressors. The results of this physical CO$_2$ flow simulation allow to assess the technical feasibility based on parameters such as velocity or pressure.
A major challenge in using such a detailed planning tool for a very long-term perspective is the definition of input data, most specifically the future CO$_2$ quantities at very high spatial resolution. In addition, the specific planning of pipelines requires developing a general topology of an overall pipeline infrastructure system that combines the individual sources and sinks of CO$_2$. This draft topology can be used as a starting point for the technical planning of the pipeline system with the software MYNTS. This method can be applied generally as long as the described data is available in sufficient resolution. In this paper, we will apply this methodology to a German case study, the details of which are described below.
In the following, we describe the four steps of our method:
\begin{enumerate}
    \item \textbf{Scenario definition:} Sets the broad frame and defines the system boundary. The scenario describes a possible future coherent CO$_2$-neutral energy system in Germany and defines components of the CO$_2$ system.
    \item \textbf{Calculation of spatial CO$_2$ balance:} CO$_2$ emissions from main sources and sinks and allocation to individual industrial sites.
    \item \textbf{Development of a CO$_2$ pipeline topology:} Connection of individual points to a German-wide pipeline system aiming at least distances and following existing gas pipeline corridors.
    \item \textbf{Physical simulation:} Describes the algorithm for the physical CO$_2$ flow modeling and the design of the specific system components (e.g. pipeline diameters).
\end{enumerate}

\subsection{Scenario description}\label{ch2_scenarios} 
As starting point for the analysis of CO$_2$ infrastructure we use a consistent industry decarbonisation scenario that was originally developed in the frame of the research project "Langfristszenarien III", calculated for the German Ministry of Economic Affairs and Climate Protection (BMWK) \cite{Fleiter.2024}. \\
The scenario originally named "O45-H2" aims to achieve climate neutrality in Germany by 2045, utilizing hydrogen-based processes as key decarbonization strategy. The scenario prioritizes the phase-out of fossil fuels by switching to hydrogen or electricity where more efficient. It applies carbon capture (CC) to hard-to-abate sectors that lack alternative abatement options. 
Specific CO$_2$ sources considered in this scenario are cement and lime production (\ref{ch2.1.1}) and waste incineration (\ref{ch2.1.2}). CO$_2$ sinks considered are Carbon Capture and Utilisation (CCU) to produce basic chemicals (olefins \& aromatics, methanol, urea) (\ref{ch2.1.3}) and offshore geological CO$_2$ storage under the North Sea (\ref{ch2.1.4}). Industry site information is derived from Neuwirth et al. \cite{Neuwirth.2022}.
As such, the scenario is in line with recent political discussions and the focus of Germany's Carbon Transport and Storage Law (KSpTG) \cite{FederalGovermentGermany.2025}. The scenario does not consider any production of e-fuels for the transport sector or the facilitation of negative emissions through direct air capture (DAC) or biomass use in Germany. In the following, we refer to this scenario as "HtA+CCU". 

Exploratory simulation runs have shown that the large-scale use of CCU substantially reduces the quantities of CO$_2$ transported, because average distances to major sinks are shorter compared to offshore storage. At the same time, the large-scale development of CCU in Germany is highly uncertain.  

 To make sure that the analyzed pipeline system will also work under such uncertainty, we develop an adapted version, the "HtA+Imports" scenario. Besides the exclusion of CCU, the adapted scenario considers two additional changes that increase the volume of transported CO$_2$ and make it a more robust solution. These relate to CO$_2$ cross-border flows and plastics recycling. The original scenario does not consider CO$_2$ imports from neighbouring countries to Germany, while the central location of Germany in Europe makes such future cross-border flows very likely.  They will further increase the volume of transported CO$_2$. Furthermore, the "HtA+Imports" scenario considers higher emissions from waste to energy (WtE) plants, due to less ambitious assumptions for progress in plastics recycling compared to the original scenario.

\subsubsection{Cement and lime production}\label{ch2.1.1}
The energy-intensive production of cement and lime generates significant CO$_2$ emissions, mainly from the calcination of limestone. About two-thirds of emissions are process-related and difficult to avoid \cite{Mohr.2024}. The current energy supply relies heavily on fossil fuels, but also contains a high share of secondary waste-based energy carriers, mostly of fossil nature \cite{BVKalk.2024, Mohr.2024}. Several measures can contribute to decarbonization of cement production, but even with ambitious assumptions, it is unlikely that near-zero CO$_2$ cement production is possible in Germany by 2045 without carbon capture. However, measures like fuel switch, material efficiency, new clinker and cement receipts and use of secondary/alternative materials have the potential to substantially reduce emissions from cement production and, thus, reduce the remaining need for carbon capture. We make the following specific assumptions about the individual abatement strategies.
\begin{itemize}
    \item \textbf{Fuel switch in process heat:} Electrification, use of hydrogen, biomass, and alternative fuels replace natural gas and coal.
    \item \textbf{Material efficiency \& substitution:} Cement production is assumed to decline by 20\% by 2045, and clinker production by 33\%, due to the use of low-clinker cements.
    \item \textbf{New binders \& substitute materials:} Clinker is partially replaced by alternatives (e.g., pozzolans, sewage sludge ash), though limited by quality standards and raw material availability.
    \item \textbf{Decline in lime production:} A 40\% reduction by 2045 due to lower demand in the energy and steel sectors.
    \item \textbf{CCUS:} CO$_2$ capture and storage or utilization (CCS/CCU) as a complementary strategy to achieve climate neutrality. We assume carbon capture to be deployed at all cement and most lime plants latest by 2045.
\end{itemize}

\begin{figure}[H]
        \centering
        \includegraphics[width=0.7\linewidth]{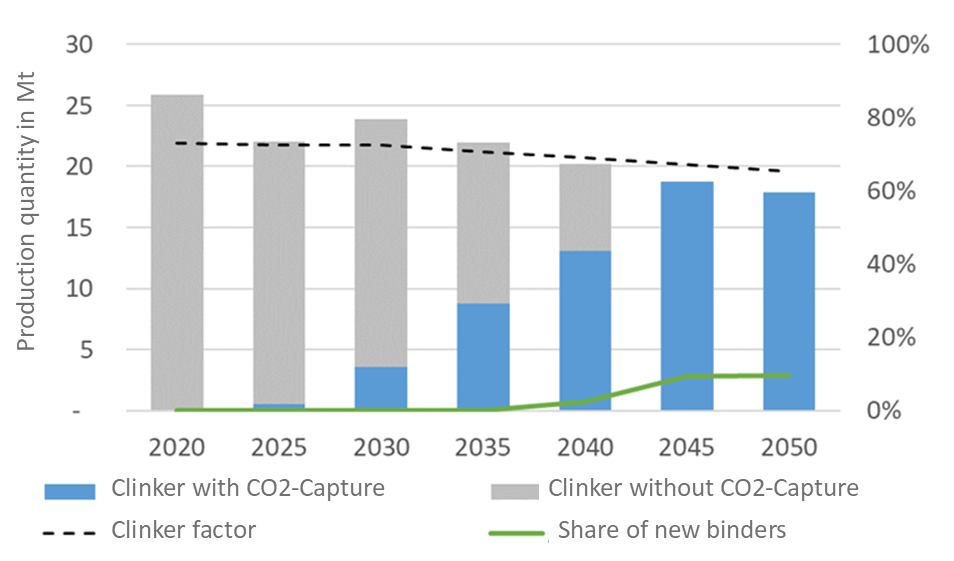}
        \caption{Development of cement clinker production in the scenario (adapted from \cite{Fleiter.2024})}
        \label{fig:cement_development}
\end{figure}

Figure \ref{fig:cement_development} shows the assumed uptake of carbon capture at cement plants, alternative binders, and the expected clinker production quantities. Both scenarios use the same assumptions for captured emissions from cement and lime.

\subsubsection{Waste incineration}\label{ch2.1.2}
Waste incineration emits CO$_2$ coming from both renewable and fossil sources. An examples for fossil sources are most plastic packaging or plastic products, whereas renewable sources can come from food waste or biobased plastics products. Mitigating the fossil emissions from waste incineration will require substantial changes in e.g., the plastics value chains and raw materials.
In the original "HtA+CCU" scenario emissions from plastic waste incineration are considered to decrease due to increased recycling progress. 
In the adapted "HtA+Imports" scenario we assume that this progress will not take place and emissions from waste incineration plants remain on the level of the year 2020. We further assume all major waste incineration plants to use carbon capture with a capture rate of 90\% in both scenarios.

\subsubsection{Basic chemical production} \label{ch2.1.3}
High-value chemicals (HVCs), such as ethylene, are crucial building blocks for the production of plastics and various everyday products. Currently, HVCs are primarily produced through the steam cracking of fossil naphtha, a byproduct from refineries, with ethylene serving as a representative for all HVCs. To achieve climate neutrality, several alternative production routes and mitigation strategies are being explored in the scenarios as outlined in Table 
\ref{tab:chemical_industry_scenarios}.
 Due to the missing regulatory incentives to defossilize feedstocks, the high abatement costs and the resulting uncertainty about the future uptake of CCU, we exclude it in the finally calculated scenario \cite{Neuwirth.2024}.

\begin{table}[H]
    \centering
    \caption{Strategies for reducing emissions from the production of high-value chemicals}
    \label{tab:chemical_industry_scenarios}
    \begin{tabular}{p{3.5cm}p{4.5cm}p{4.5cm}}
        \hline
        & \textbf{Scenario HtA+CCU}& \textbf{Scenario HtA+Imports} \\

        \hline
        Decarbonizing Process Heat & Electrification of the steam cracker. & Electrification of the steam cracker. \\
        \hline
        Material Efficiency and Substitution & Decline of 8\% in plastics production by 2045 & Decline of 8\% in plastics production by 2045 \\
        \hline
        Mechanical Recycling & Increase in mechanical recycling, by 6\% & No change compared to 2020 recycling rates\\
        \hline
        Chemical Recycling & No chemical recycling & No chemical recycling \\
        \hline
        CCU & Synthesis of HVCs from methanol via methanol to olefins (MTO) and methanol to aromatics (MTA), using methanol produced from captured CO$_2$ and hydrogen. Import of Fischer-Tropsch fuels and naphtha & Continued use of electrified steamcrackers with fossil fuels instead of CCU and methanol based production \\
        \hline
    \end{tabular}
\end{table}




\subsubsection{\COtwo-Hubs}\label{ch2.1.4}
For connection the CO$_2$ network to storage facilities  four ports (Wilhelmshaven, Bremerhaven, Brunsbüttel, and Rostock) are being considered \cite{AFRedaktour.22.03.2022, BenjaminKlare.07.02.2024, Equinor.04.07.2022, wintershalldeaHESWilhelmshaven.18.10.2022}. These ports intend to install CO$_2$ terminals for loading ships, which will then head to offshore storage facilities in the North Sea. Furthermore, import and export of CO$_2$ from and to neighboring countries must be considered. For this purpose, border crossing points have been assumed at the location of the current gas border crossing points. This allows mapping of CO$_2$ flows through Central Europe and the shared use of storage capacities planned in the Netherlands, Norway, and Denmark.

\subsection{Resulting \COtwo balance}

\begin{figure}[H]
        \centering
        \includegraphics[width=1\linewidth]{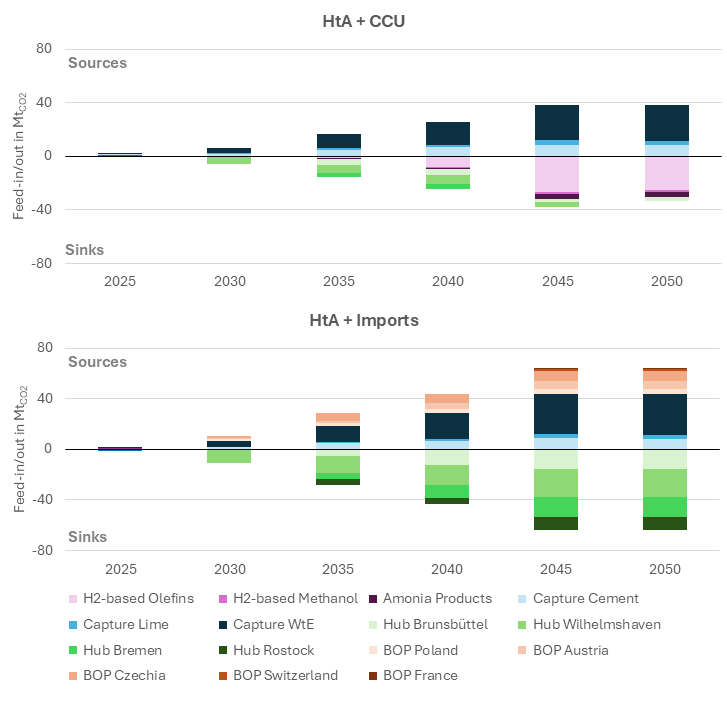}
        \caption{CO$_2$ balance and temporal development of both scenarios}
        \label{fig:temp_dev}
\end{figure}

The resulting CO$_2$ balance including major source and sink categories is shown as temporal development in Figure~\ref{fig:temp_dev} in comparison of the two scenarios. The scenario "HtA+CCU" shows an annual CO$_2$ demand of 31 Mt in 2045, with quantities of 27 Mt CO$_2$ required for the production of methanol alone. This demand reflects a rapidly increasing production volume of green methanol. Large parts of the captured emissions from lime, cement  and waste incineration plants are used for the production of green chemical feedstocks leaving only little need for CO$_2$ storage.
The  scenario "HtA+Imports", which does not envisage the use of CO$_2$ in the chemical industry, shows significantly larger quantities that need to be transported to storage sites. In addition to cement, lime, and waste incineration emissions (8.5, 3.25, and 32 Mt/a respectively in 2045), this also includes transit flows from neighbouring countries amounting to 20 Mt/a in 2045. These emissions are sent to storage sites under the North Sea via four CO$_2$ hubs.

The quantities from 2045 in scenario "HtA+Imports" are used as the basis for the following calculations. A network that is sized to transport the quantities in this scenario will likely also work for the scenario "HtA+CCU", with its smaller transport quantities. 

Table \ref{tab:scenarios} provides an overview of the quantities that need to be transported in the two scenarios in 2045.

\begin{table}[H]
    \centering
    \caption{CO2 quantities for individual sources and sinks for both scenarios in 2045 [in Mt/a]}
    \label{tab:scenarios}
    \begin{tabular}{cccccccc}
        \toprule
        Scenario & Name & Import & Export & HtA & CCU & WtE & Hubs\\
        \midrule
        1 & HtA+CCU& - & - & 11.3 & 31.92 & 26.7 & 6.71 \\
        2 & HtA+Imports & 20 & - & 11.3 & - & 32.04 & 63.34 \\

        \bottomrule
    \end{tabular}
\end{table}


\subsection{Development of \COtwo pipeline topology}
\label{pipeline-topology}
Based on the definition of potential CO$_2$ sources and sinks for both scenarios we derive a potential CO$_2$ pipeline topology. 
Due to the difference in density, dense-phase CO$_2$ transport allows much higher transport volumes than gas-phase CO$_2$ transport, which makes it the preferred solution for long-distance transport pipelines. Dense-phase transport, however, is not possible with existing natural gas pipelines because the wall thickness cannot withstand the significantly higher pressures \cite{AIT2024_CO2Netz}. Thus, we assume that the entire network will consist of newly constructed dense-phase CO$_2$ transport pipelines

Converting existing natural g as pipelines and using them for gaseous transport could lead to cost reductions, but would require more specific assessment beyond the scope of our analysis  \cite{Schützenhofer.10.2024}. 
Still, we assume the pipelines to follow the routes of the current natural gas transmission network to reduce the costs and time required for approval and construction. 
 We apply a two step heuristic approach: 
\begin{enumerate}
    \item \textbf{Connect each site} to the nearest natural gas transmission network pipeline corridor directly.
    \item \textbf{Connect sinks and sources} by following the routes of the natural gas transmission network pipelines.
\end{enumerate}

\begin{figure}[!ht]
        \centering
        \includegraphics[width=1\linewidth]{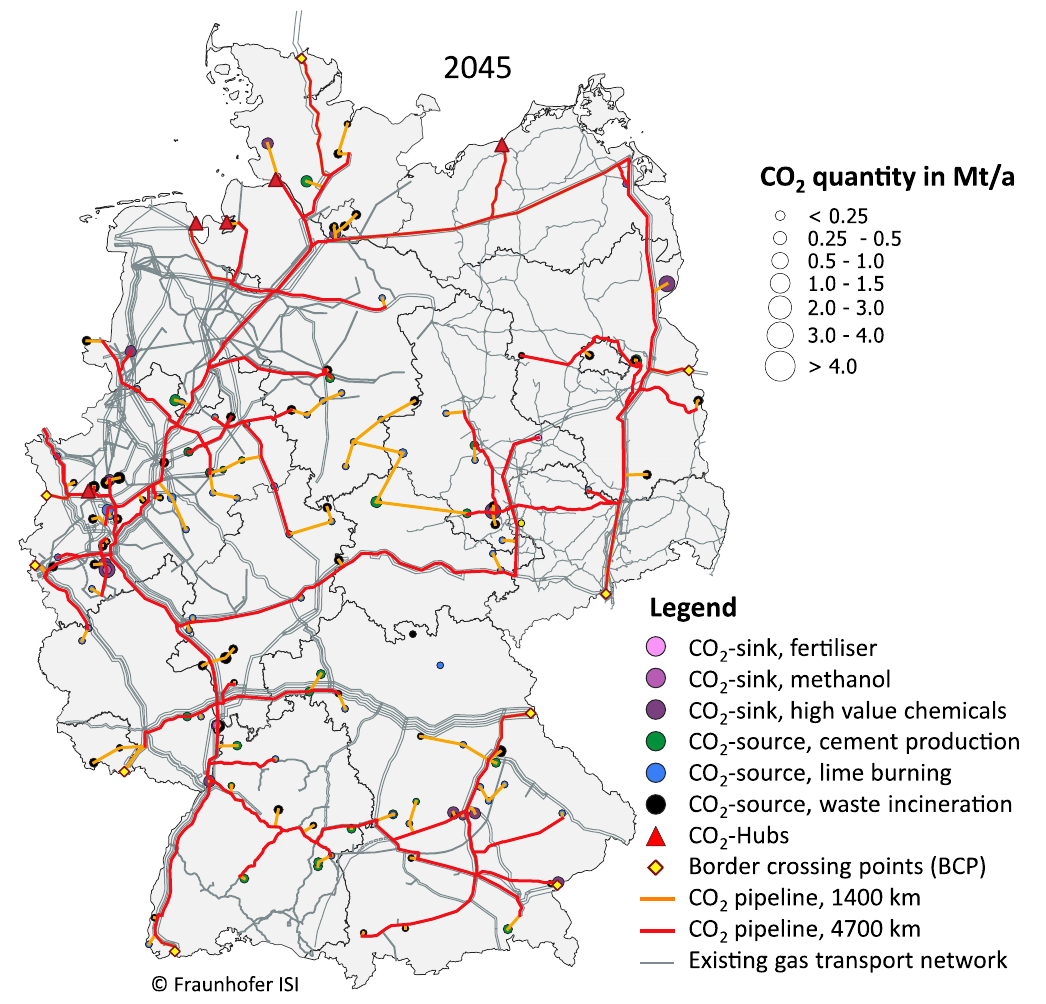}
        \caption{CO$_2$ Topology (adapted with modifications from \cite{Fleiter.2024})}
        \label{fig:co2_topology}
\end{figure}
Figure~\ref{fig:co2_topology} shows the resulting pipeline topology. An earlier version of this network was already used by Fleiter et al. as part of the “Langfristszenarien III” project for the German Ministry of Economic Affairs and Climate Protection (BMWK) \cite{Fleiter.2024}.
The proposed CO$_2$ gird should not be understood as an analysis or projection for individual sites, but rather be a rough first estimation on a potential network. 
The CO$_2$ sinks reflect the current refinery and chemical sites as well as CO$_2$ hubs for sequestration. The current locations of lime (52 sites),
cement (32 sites), and waste incineration (55 sites) are taken into account. We do not consider any potential spatial shifts, new establishments, or closures. 
Based on the above assumptions, the pipeline length is approximately 7,000 km. The transport is assumed to take place completely in the dense phase. The necessary liquefaction takes place at each feed-in site.

In order to simplify the simulation task generated from such a network, we performed the following steps: For each part between two junctions $A$ and $B$ where more than two pipelines meet, assume a straight pipeline, whose length is the sum of all parts between $A$ and $B$ in the original gas network. The heights of $A$ and $B$ can be left unchanged. For a real pipeline planning task, all points between with significant height deviations should have been included, which was not necessary for this proof-of-concept study. This is the main reason why the network graphs in the following (from the simulation tool) look somewhat different compared to Figure \ref{fig:co2_topology}.

\subsection{Simulation and optimization of \COtwo infrastructure network}
The potential future CO$_2$ pipeline topology was introduced in Sec.~\ref{pipeline-topology} and is based on the CO$_2$ sources and sinks of the scenario "HtA+Imports". However, the feasibility of this network, as well as the physical characteristics of CO$_2$ infrastructure necessary for successful dense phase transport — without phase changes during transit that could lead to cavitation and damage to the pipeline \cite{raimondi2022ccs, sleiti2022co2} — has not yet been explored. Therefore, we use a physical pipeline simulator that allows us to specify the technical dimensioning of the system and test its feasibility. We do this by addressing the following questions:
\begin{itemize}
    \item What are the technical parameters (or physical criteria) — such as pipeline diameter and thickness, temperature and pressure conditions, the number of required pumps, and overall infrastructure design — for the CO$_2$ pipeline topology introduced in Sec.~\ref{pipeline-topology} that enable the transport of CO$_2$ in dense phase?
    \item What are the technical characteristics of the CO$_2$ infrastructure network for cost-effective investment?
\end{itemize}

There are several studies simulating the transport of CO$_2$ through pipelines and its pumping to underground storage, taking into account the potential phase transition, specifically the conversion from dense to gas phase due to the presence of impurities in the network \cite{aursand2013pipeline, chaczykowski2012dynamic, liljemark2011dynamic, nimtz2010modelling}.  Most of these studies focus on transport within a single pipeline rather than addressing larger networks, such as a national CO$_2$ infrastructure network. 

To tackle the aforementioned questions and simulate the German national CO$_2$ network, we have developed the software MYNTS (Multiphysical Network Simulator) \cite{baldin2021advanced,baldin2020topological,clees2016making, clees2016mynts}, which is designed to plan, simulate, and optimize energy networks across various sectors, including gas, electricity, heat, and water.

In the following subsections, we will explain the software MYNTS and discuss the chosen and required parameters for simulating the CO$_2$ network as well as demonstrating its feasibility.

\subsubsection{Network simulation software: MYNTS}
For the simulation and analysis of the German CO$_2$ infrastructure network, the multi-physics software MYNTS has been enhanced to model CO$_2$ transport in dense phase, as well as phase transitions due to impurities or the mixing of CO$_2$ with other components in pipelines.

For stationary fluid transportation problems involving both liquids and gasses, we implement the following standard pipe transport equations.
\begin{equation}
    \frac{dP}{dx}=\frac{-\lambda}{2D}\rho v^2-\rho v\frac{dv}{dx}-\rho g\frac{dh}{dx}
    \label{pressure}
\end{equation}
for pressure $P$, where $x$ is the dimension of space in the direction of the pipe, $\rho$ and $v$ are the density and velocity of the fluid, $D$ is the inner diameter of the pipe, $g$ is the gravitational acceleration, and 
$h$ is the height. On the right-hand side, the first term is usually dominant (mainly in the gas phase), representing the contribution of the friction force, defined in terms of the dimensionless friction coefficient $\lambda(k/D,Re)$ using the Nikuradse \cite{nikuradse1950laws} formula or the more accurate Hofer \cite{hofer1973error} formula. Here, $k$ is the pipe's
roughness, $Re=4|Q_m|/(\pi \mu_{visc}D)$ is the Reynolds number,
where $\mu_{visc}$ is the dynamic viscosity, and $Q_m=\rho v\pi D^2/4$ is
the mass flow, assumed to be constant along the pipe. However, in the dense phase the influence of the rightmost term becomes dominant, especially when the pipes must cross hills and mountains. This will be discussed in subsection \ref{sec:elevation}.



The temperature profiles are represented by the following equation:
\begin{equation}
    {Q_m}\frac{dH}{dx}=-\pi Dc_h(T-T_s)
    \label{temperature}
\end{equation}
where $H$ [J/kg] is the specific enthalpy, defined as $H=e+P/\rho$. In this equation, $e$ represents the specific energy, $c_h$ is the heat transfer coefficient, $T_s$ is the soil temperature. Note that when the heat exchange is turned off ($c_h=0$), for example for regulators, the process described by this formula is isoenthalpic ($dH=0$), and the temperature change is related to the pressure change by the well-known equation 
$dT=\mu_{JT}dP$, where $\mu_{JT}=(-\partial H/\partial P)_T/(\partial H/\partial T)_P$ is the Joule-Thomson coefficient.

In addition to eq.~(\ref{pressure}) and (\ref{temperature}), we utilize the Equation of State (EoS), a so called gas law, specifically the GERG2008 approximation \cite{kunz2012gerg}, in our calculations which also works well for mixtures with the typical impurities.
Kirchhoff’s laws are used for connecting the components of the
network (pipes, pumps, regulators, valves, heaters, coolers, etc.).
For more details on modeling gas transport networks and their components, see \cite{domschke2017modellierung} for a description of all relevant elements, including
regulators, pumps, heaters and coolers, and \cite{benner2018gas} and the references given therein for an index-optimal spatial discretization technique for equations~(\ref{pressure}) and (\ref{temperature}), as well as the equation of state (EoS).

By utilizing the MYNTS software, one can solve such non-isothermal 1D steady-state pipeline flow  equations. MYNTS employs a semi-coupled approach, iterating these equations alongside the equations for enthalpy and gas composition. The EoS is incorporated into this process. Using the current pressure ($P$), temperature ($T$), and gas composition represented by the vector of mole fractions ($\zeta$), the software first computes the molar density ($\rho_m$), followed by the compressibility factor ($z$) and height ($h$).
In particular, a theoretical foundation was established \cite{clees2021efficient} that ensures the existence and uniqueness of solutions in practically relevant scenarios, along with a corresponding numerical method for constructing these solutions. To accelerate simulations, topological reduction methods were developed \cite{baldin2019topological}, and model-order techniques are explored in \cite{grundel2014model}.

With the specific aim of modeling cross-sectoral applications, a so-called universal translator (UT) \cite{baldin2018universal} was developed. The UT takes a description of the network along with its scenario data (such as pump settings and profiles) and a translation matrix (TM) to formulate the system to be solved. The TM includes formulas presented in a readable format. Currently, the output can be provided for a nonlinear programming solver (e.g., IPOPT) as well as in Mathematica or MATLAB format. A proprietary nonlinear solver is utilized.

In the following sections, we will discuss the thermodynamic models and the parameters selected to simulate the CO$_2$ infrastructure network.

\paragraph{Gas-law and friction model}
As mentioned earlier, we utilize the GERG-2008 gas law \cite{kunz2012gerg}, which provides a significantly more complex description compared to other models, such as Papay \cite{papay1968atermelestechnologiai} and DC92 \cite{iso12213-2_2006}. 
In the GERG model, the compressibility factor $z$ is derived from the Helmholtz free energy, which is divided into components for the ideal gas and the residual fluid. A substantial amount of measurement data from various gas combinations is employed to fit the numerous coefficients. For further details, refer to \cite{kunz2012gerg} and ISO 20765–2/3. Importantly, GERG-DLL allows us to calculate phase transitions that can occur in the pipeline during the transport of non-pure CO$_2$ in the dense phase.

In this study, we use the Hofer formula for calculating the friction coefficient in Eq.~(\ref{pressure}). Compared to the simpler Nikuradse formulation, the Hofer approximation provides higher accuracy across a wider range of flow conditions while remaining computationally efficient \cite{benner2018gas}. For high Reynolds numbers, both formulations yield nearly identical results.

\paragraph{Pump model}
In MYNTS, pumps/compressors and their drives are modeled according to their characteristic curves describing the dependencies of flow, rotational speed, and compression \cite{clees2018modeling}. Usually, an individual set of characteristics is supplied for each device. However, since the CO$_2$ infrastructure network is a novel subject and there is currently no such network in Germany, we have not yet been able to identify the typical characteristics for the designed CO$_2$ pumps. In this study, we employ a ``free" pump model, in which the pumps are not constrained by characteristics but are assumed to raise the outlet pressure directly to the prescribed setpoint. This modeling choice enables a systematic analysis of the network’s hydraulic behavior without the limitations imposed by device-specific performance curves. While this abstraction neglects efficiency considerations and maximum compression ratios, it provides a suitable basis for exploring the fundamental operational feasibility of CO$_2$ transport within the designed infrastructure.

\paragraph{Physical parameters for the CO$_2$ network}
We use the CO$_2$ balance and locations of the scenario ``HtA+Imports" as a basis, which is defined to represent a case with a relatively high transport volume. As such, the network will be flexible enough to accommodate other scenarios and developments. As a reminder, in this scenario, 40 Mt/a of CO$_2$ will be transported through the network. Of this amount, 20 Mt/a will be imported to Germany, with no exports planned.

Table \ref{physical_parameters} presents the assumptions and variations of the physical parameters, including internal diameters, pipeline roughness, nodal heights, and other relevant quantities. Since the network does not represent an existing topology, it has been modeled in a relatively coarse resolution with 326 pipe elements. The runtime with MYNTS 4.9.5 was about 10 seconds for pure CO$_2$.   

\begin{table}[h!]
\caption{Physical parameters of the system}
\centering
\begin{tabular}{p{2.9cm} p{2.8cm} p{2.5cm} p{2.25cm}}
\toprule
Parameter               & Symbol [unit]             & Value & Variations \\
\midrule
Total Pipe Length       & $L_{\text{total}}\; [km]$   & 6937.6            & --- \\
Pipe Roughness          & $k\; [mm]$                  & 0.036             & 0.01, 0.05, 0.1 \\
Pipe Diameter           & $D\; [mm]$                   & 378 - 661     & 468, 661, 882 \\
Nodal Height         & $h\; [m]$                      & terrain height & -1 ... 1015 m\\
Heat Transfer\newline Coefficient & $c_h\; [W/(m^2 \cdot K)]$      & 2            & --- \\
Soil Temperature        & $T_s\; [^\circ C]$              & 10                & --- \\
Molar Mass              & $\mu\; [kg/mol]$            & 0.044             & --- \\
Dynamic Viscosity       & $\mu_{\text{visc}}\; [Pa \cdot s] $     & 0.0001            & --- \\
Fluid Composition       & $(\text{CO}_2,\ \text{H}_2,\ \text{O}_2)\; [\%]$ & (100, 0, 0) & (98, 1, 1),\newline (96, 2, 2)\\
Setpoint Pressure         & $P_{\text{set}}\; [barg]$   & 125              & --- \\
\bottomrule
\end{tabular}
\label{physical_parameters}
\end{table}

\section{Results}
In this section, we present our simulation results and discuss the network parameters that enable the feasible transport of CO$_2$ in dense phase, as well as their impacts on investment costs. All simulations are based on the scenario "HtA+Imports".


\subsection{Network Costs according to Pipeline Diameter on Pump Requirements in \COtwo Transmission Systems}

The network design initially assumed uniform pipeline diameters throughout. Three standard diameters were evaluated: DN500, DN700, and DN900. While total pipeline costs increase with diameter, the maintenance costs of the pumps decrease, since larger pipelines require fewer pumps in the system. In a follow-up consultation with a German transmission system operator, it was noted that the use of DN900 pipelines for dense phase CO$_2$ transport is generally not considered realistic under current industry practice. Nevertheless, DN900 is retained in the analysis to assess the theoretical reduction in the number of pumps that such large diameters could enable. Based on informed estimates and preliminary cost assessments, the optimal configuration was identified as using DN700 pipelines in sections with high mass flow and DN500 or DN400 pipelines in areas with significantly lower flow, typically near the ends of branches. The total costs for this configuration are comparable to those of a network using only DN500 pipelines but offer the advantage of requiring fewer pumps (see Table \ref{pipeline_costs}).

Specific pipeline cost estimates were derived from the Network Development Plan 2022 by FNB Gas \cite{fnbgas2022}, which provides data for newly constructed hydrogen pipelines. These values were increased by approximately 50 \% to account for the higher material costs expected for CO$_2$ pipelines, primarily due to thicker wall requirements and the use of higher-grade steels. Additionally, according to the CO$_2$ Network Feasibility Study for Austria \cite{AIT2024_CO2Netz}, the specific investment costs are given as 2115 €/m for DN400, 2355 €/m for DN500, and 2955 €/m for DN700 (no data reported for DN900). These values align well with the estimated cost adjustments and suggest that both sources produce similar specific cost levels for CO$_2$ pipeline infrastructure.

By analyzing various diameter configurations, we examined how pipeline sizing affects the number of required pumps and the overall network cost. Smaller diameters may lower material costs but lead to higher operational and maintenance expenses due to increased pumping demands. Conversely, larger diameters can reduce these operational requirements at the expense of higher investment costs.

According to the CO$_2$ Network Feasibility Study for Austria \cite{AIT2024_CO2Netz}, the investment costs for individual pumps vary depending on the mass flow, manufacturer, and the number of pumping units installed at each station, ranging between approximately 3.6 million € and 13.2 million €. These costs account for the total expenses associated with a complete pump station, including auxiliary equipment, control systems, and installation. To validate these estimates, a German transmission system operator was consulted, who reported a representative average pump station cost of about 7.5 million €. Consequently, this value was adopted as the reference for the calculations presented in Table \ref{pipeline_costs}.

\begin{table}[h!]
\caption{Network topology related costs}
\centering
\begin{tabular}{p{3.5cm} p{1.5cm} p{2.5cm} p{1.5cm} p{1.5cm}}
\toprule
Pipeline DN             & DN500 & DN700 \&\newline (DN500/DN400) & DN700 & DN900$^{*}$ \\
\midrule
Inner Diameter $[mm]$     & 468   & 661 \&\newline (468/378) & 661 & 882 \\
Max Velocity $[m/s]$      & 6.1   & 4.0 & 4.0 & 2.1 \\
Amount of Pumps           & $>20$ & 9 & 9 & 0 \\
Specific Cost of Pumps $[\text{mil. € per Station}]$ & 7.5 & 7.5 & 7.5 & 7.5\\
Total Cost of Pumps $[\text{mil. €}]$ & 150 & 67.5 & 67.5 & 0\\
Specific Cost of Pipelines $[\text{€}/m]$ & 2500 & 3000 \&\newline (2500/2200) & 3000 & 3500 \\
Total Cost of Pipelines $[\text{bil. €}]$ & 17.3 & 17.3 & 20.8 & 24.3 \\
\bottomrule
\end{tabular}

\begin{flushleft}
\footnotesize $^{*}$DN900 pipelines for dense phase CO$_2$ transport are generally not considered realistic according to consultation with a German transmission system operator; the case is included only for theoretical comparison and to provide a broader parameter space for optimization.
\end{flushleft}

\label{pipeline_costs}
\end{table}

\begin{figure}[H]
    \centering
    \begin{subfigure}[t]{0.48\linewidth}
        \centering
        \includegraphics[width=\linewidth]{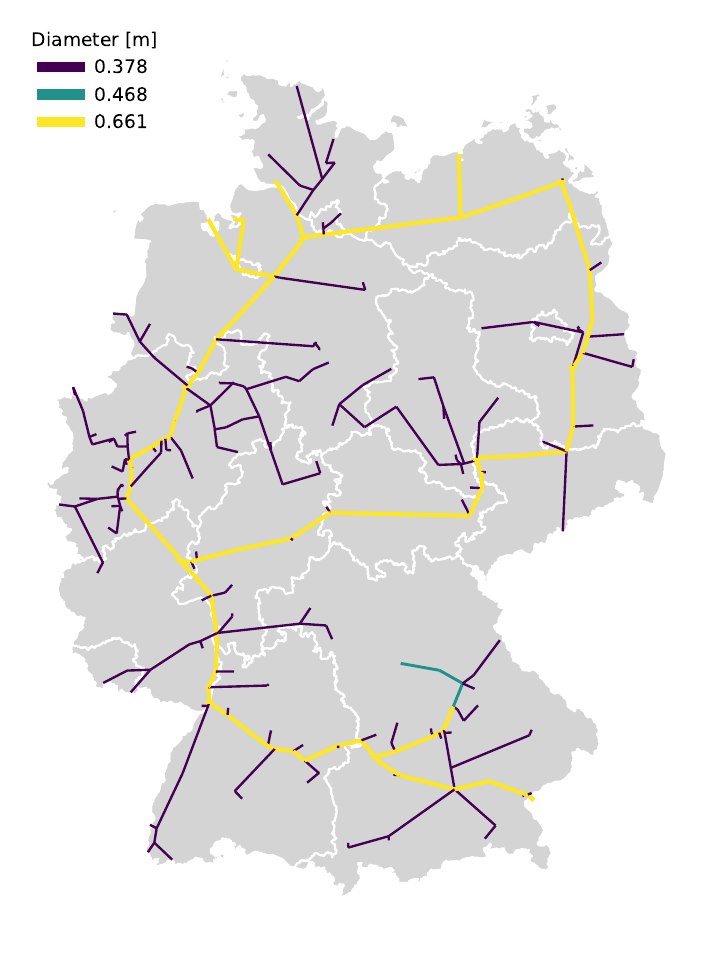}
        \caption{Pipeline Diameters}
        \label{fig:topo_dia}
    \end{subfigure}
    \hfill
    \begin{subfigure}[t]{0.48\linewidth}
        \centering
        \includegraphics[width=\linewidth]{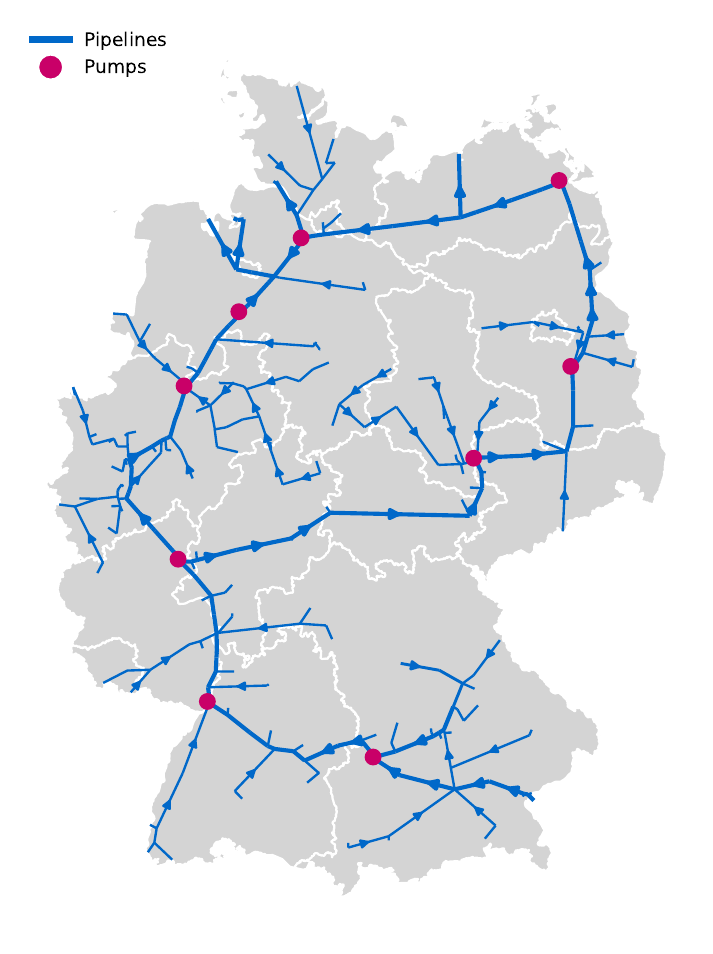}
        \caption{Pump Locations}
        \label{fig:topo_pumps}
    \end{subfigure}
\caption{Network topology showing (a) pipeline diameters and (b) pump placements, optimized for cost and hydraulic performance. Arrows indicate the flow direction.}
\label{fig:topology_pics}
\end{figure}

Figure \ref{fig:topology_pics}(a) illustrates the resulting distribution of pipeline diameters within the optimized network. Larger diameters (DN700) are concentrated along the main pipelines, where mass flow rates are highest and pressure losses would otherwise be significant. Toward the outer branches of the network, smaller diameters (DN500 and DN400) are used to reduce material costs while maintaining acceptable flow velocities and pressure levels. This configuration achieves an effective balance between hydraulic performance and overall system cost.

To locate pumps in the CO$_2$ network, we first analyzed the flow direction and determined the expected distribution across the different branches of the system. The overall flow originates predominantly in southern Germany and divides near the border between Rhineland-Palatinate and Hesse, where it continues northward through two distinct branches that eventually merge again in the Hamburg region. This flow structure provided a natural basis for dividing the network into several sections.

For each section, we evaluated how far dense phase CO$_2$ could be transported within the predefined pressure limits (minimum 85~barg, maximum 180~barg). Based on this method and the corresponding pressure behavior, pumps were placed at locations where pressure recovery was required. Since maintenance costs increase with the number of pumps installed, the network was optimized to operate with as few pumping stations as possible while still ensuring reliable transport.

This approach follows the segmentation method, originally developed in electrical engineering for the analysis of microwave planar circuits by dividing complex networks into smaller, more manageable segments with simpler geometries \cite{segment_method}. Although conceived for circuit design, the same principle can be effectively applied to gas and fluid network modeling. By partitioning a pipeline network into hydraulically consistent segments, one can efficiently assess pressure losses, flow limits, and optimal equipment placement under operating constraints.

The resulting flow directions are illustrated by arrows along the pipelines, showing the CO$_2$ movement from south to north. Figure~\ref{fig:topology_pics}(b) depicts the final pump locations within the network. Pumps are positioned primarily along the main transport routes, especially after long pipeline stretches or at major elevation changes, where pressure boosting is necessary to maintain the minimum operating pressure of 85~barg. This configuration ensures continuous CO$_2$ transport throughout the network while minimizing the number of pumping stations and associated maintenance costs.

\subsection{Pressure and Temperature Distribution}

\begin{figure}[H]
    \centering
    \begin{subfigure}[t]{0.48\linewidth}
        \centering
        \includegraphics[width=\linewidth]{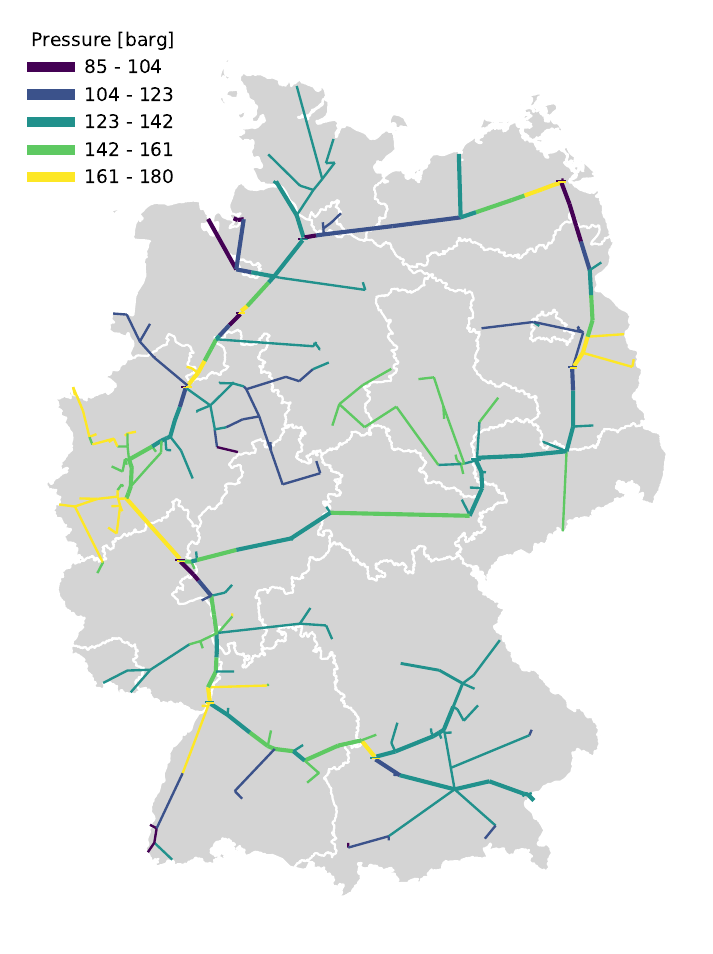}
        \caption{Pressure Distribution}
        \label{fig:res_pres}
    \end{subfigure}
    \hfill
    \begin{subfigure}[t]{0.48\linewidth}
        \centering
        \includegraphics[width=\linewidth]{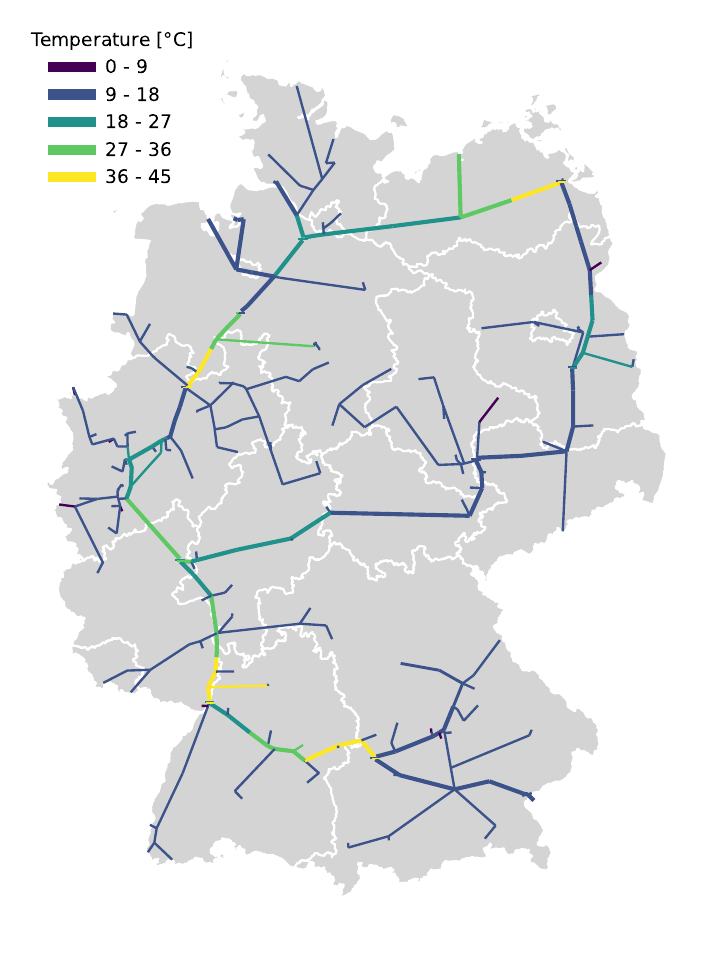}
        \caption{Temperature Distribution}
        \label{fig:res_temp}
    \end{subfigure}
\caption{Fluid dynamic simulation results of (a) pressure and (b) temperature within the optimized CO$_2$ transport network, which can be seen in Figure \ref{fig:topology_pics}.}
\label{fig:result_pics}
\end{figure}

Figure \ref{fig:result_pics}(a) illustrates the pressure distribution along the network, highlighting how pressure gradually decreases along each transport route and recovers at pump stations. The highest pressures occur immediately downstream of the pumps, where compression takes place, while the lowest pressures are found just upstream of the pumps, at the inlet points where the pressure approaches the minimum design limit of 85 barg.

To ensure that the network operates within the designed pressure and temperature limits, the temperature of the fluid must be reduced after it has been compressed. In our model, the cooling process that follows compression is represented as an output decision variable. However, in practice, cooling is associated with additional costs, which are currently unknown. This uncertainty makes it difficult to optimize the number of pumps in operation and the degree of cooling required after each pump. Moreover, the number of pumps required is directly influenced by the cooling temperature: Stronger cooling reduces the need for pumps, whereas weaker cooling increases it.

Figure~\ref{fig:result_pics}(b) shows the corresponding temperature distribution along the network. After each pump, the CO$_2$ temperature is actively controlled to represent post-compression cooling. In our simulations, we initially set the outlet temperature to 10 °C after each pump to ensure that the CO$_2$ remained in the dense phase. To explore potential cost reductions, the cooling temperature was then gradually increased, one pump at a time, up to a maximum of 45 °C where operating margins allowed. In some sections, however, such an increase was not feasible, and the temperature could only be raised to the maximum level permitted by thermal and pressure constraints.

\begin{figure}[H]
    \centering
    \includegraphics[width=\linewidth]{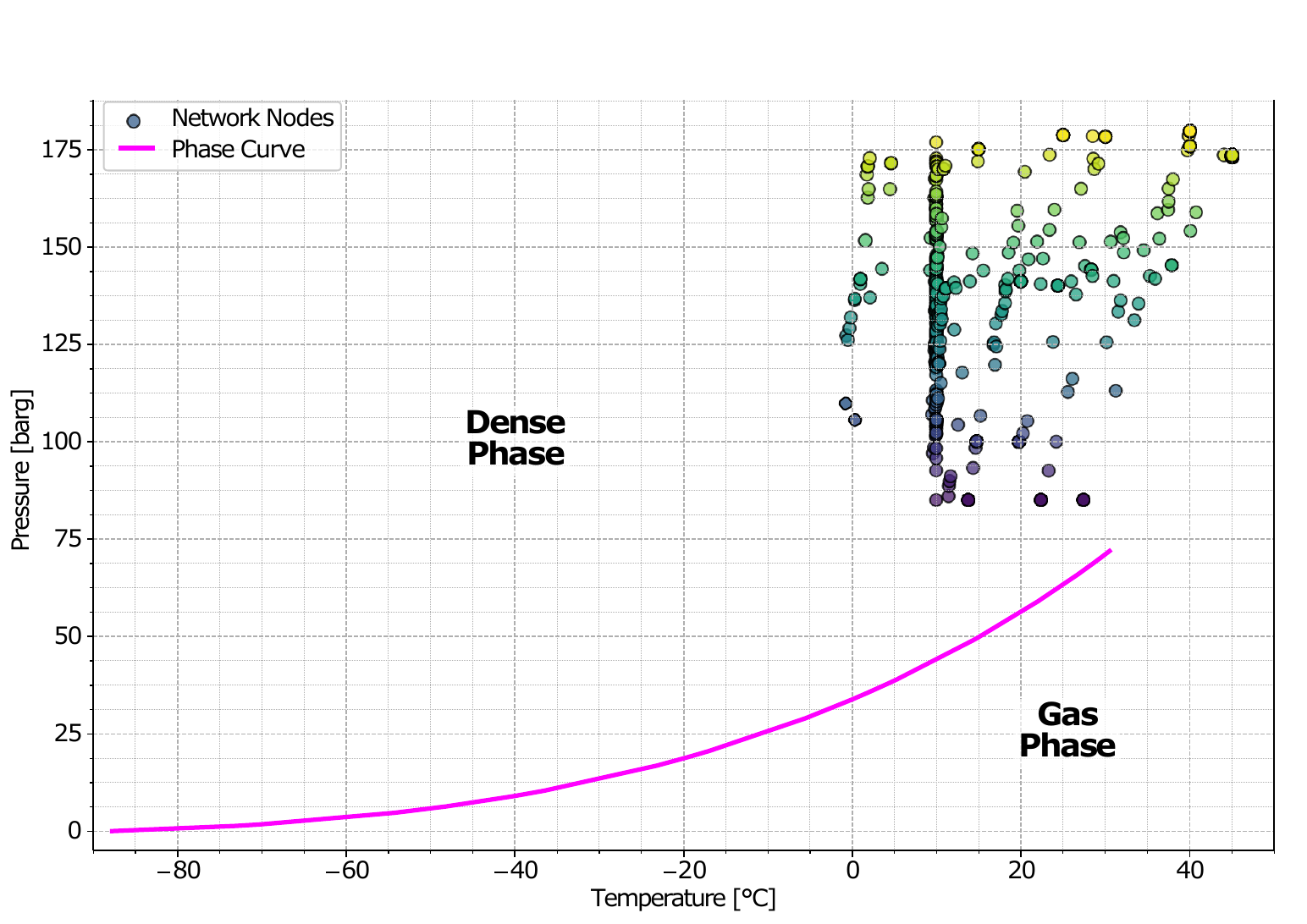}
    \caption{Phase diagram of CO$_2$ showing the pressure–temperature distribution of all network nodes. The pink line marks the boundary between the gaseous and dense phases. Each data point represents a node in the network.}
    \label{fig:phase_diagramm}
\end{figure}

Figure~\ref{fig:phase_diagramm} shows the phase diagram of CO$_2$, where the pink curve represents the boundary between the gaseous and dense phases. All network operating points must be above this line to ensure that CO$_2$ remains in the dense phase throughout the system as crossing below this boundary would lead to gas formation and potential instability in the network.

The x-axis indicates the temperature, while the y-axis shows the corresponding pressure. Each point in the plot corresponds to a network node. The color gradient of the data points reflects the pressure distribution, consistent with the results shown in Figure~\ref{fig:result_pics}(a). This visualization helps to verify that the selected cooling strategy and the number of pumps are sufficient to keep CO$_2$ within the desired thermodynamic region.




\subsection{Impact of Elevation Differences on Pressure Distribution in \COtwo Transmission Systems}
\label{sec:elevation}
Due to gravitational effects in fluid transport, we investigate how elevation differences influence pressure behavior within the network.

\begin{figure}[H]
    \centering
    \begin{subfigure}[t]{0.48\linewidth}
        \centering
        \includegraphics[width=\linewidth]{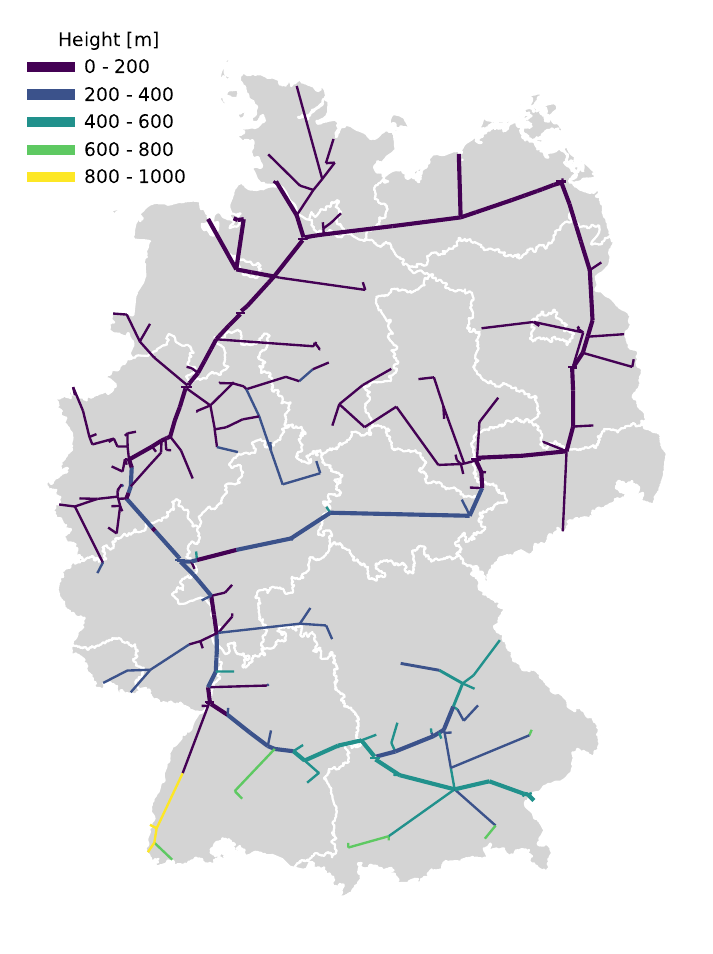}
        \caption{Nodal Elevation Distribution}
        \label{fig:height_ger}
    \end{subfigure}
    \hfill
    \begin{subfigure}[t]{0.48\linewidth}
        \centering
        \includegraphics[width=\linewidth]{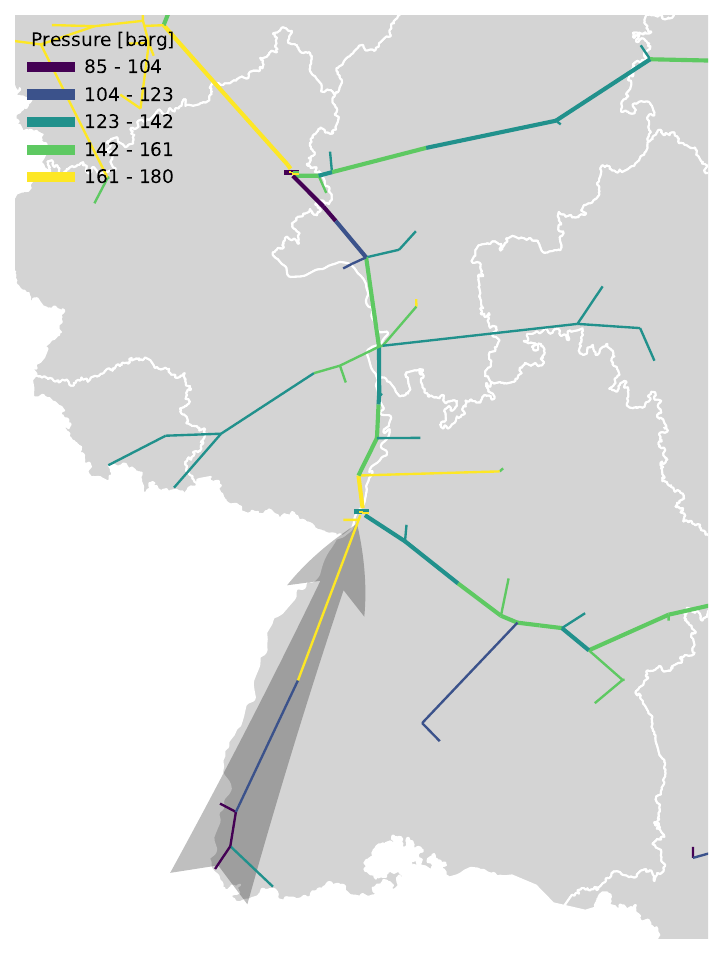}
        \caption{Upstream Lower Pressure}
        \label{fig:upstream_low}
    \end{subfigure}
\caption{Effect of elevation differences on CO$_2$ pressure, showing (a) network node elevations and (b) pressure distribution along the highlighted section influenced by the hydrostatic pressure effect. The arrow in (b) marks the pipeline segment between the higher-elevation upstream region and the lower-elevation downstream region, where this pressure behavior occurs.}
\label{fig:height_pics}
\end{figure}

The simulations account for the elevation differences between the locations where CO$_2$ is produced and where it is transported. Figure~\ref{fig:height_pics}(a) shows the nodal elevation distribution across the entire network, covering almost all of Germany and illustrating the significant variation in terrain. Figure~\ref{fig:height_pics}(b) provides a zoomed-in view of a specific region where the hydrostatic pressure effect becomes evident.

In typical gas networks with negligible elevation changes, the pressure is generally higher upstream, and gas flows from higher to lower pressure. In contrast, for dense phase CO$_2$, the pressure distribution is strongly influenced by elevation differences. As shown in Figure~\ref{fig:height_pics}(b), the fluid at an upstream location near the German–French–Switzerland border can exhibit a lower pressure than the downstream location around Karlsruhe when flowing from higher to lower elevation. This effect occurs because the weight of the CO$_2$ column increases the pressure at lower elevations within the network. The relevant pipeline section, marked by an arrow in the figure, indicates both the flow direction and the region discussed here.

The described effects of elevation are well-known, but users of simplified transport models might not be fully aware of the strong influence of hight differences, especially in context with dense phase. 

\begin{figure}[H]
    \centering
    \includegraphics[width=\linewidth]{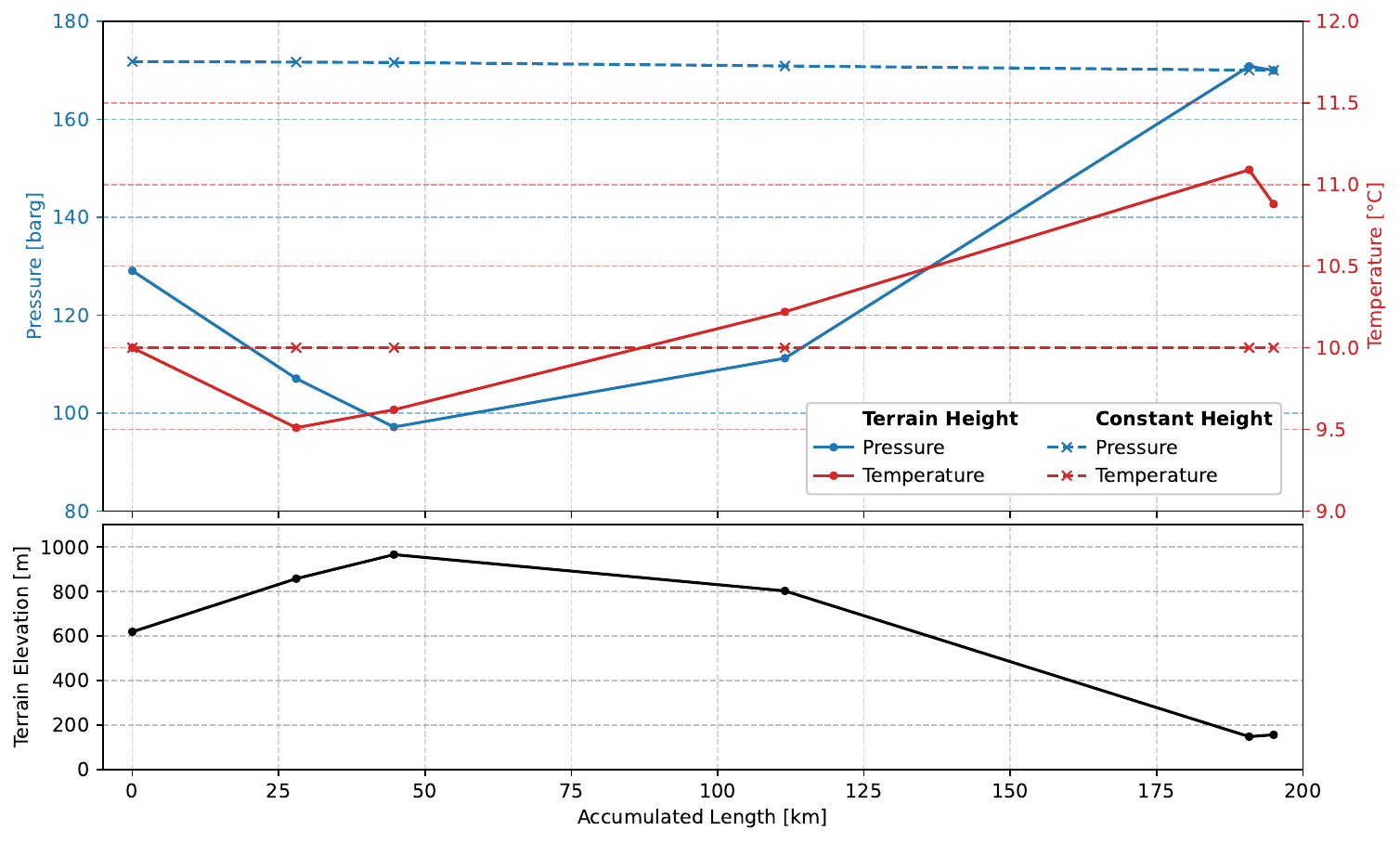}
    \caption{Comparison of pressure and temperature distribution along the selected pipeline section. Solid lines represent the simulation with actual terrain elevations, while dashed lines correspond to the constant-height (h=0) case.}
    \label{fig:no_height_custom}
\end{figure}

The Figure~\ref{fig:no_height_custom} presents the pressure and temperature behavior in greater detail for the pipeline section indicated by the arrow in Figure~\ref{fig:height_pics}(b). The solid lines represent the nodal pressure and temperature values from the simulation that accounts for actual terrain elevations, while the dashed lines correspond to the simulation in which all node elevations are uniformly set to \(h=0\).

In the terrain-height simulation, both pressure and temperature vary according to changes in elevation. In contrast, when constant node heights are assumed, the pressure drop is considerably smaller, and the temperature remains nearly constant at approximately 10 °C. The pressures in both simulations converge at the downstream end of the route (around 195 km accumulated length) because a pump located there operates under inlet control, maintaining a constant inlet pressure of 170 barg. This setting is necessary in the terrain-height case to prevent upstream pressures from falling below 85 barg. Although a lower inlet pressure setting of the pump could be used for the constant-height scenario, the same pump configuration is retained in both simulations to ensure consistent comparison and to provide a clear interpretation of the pressure–temperature behavior.

\subsection{Impact of Impurities in \COtwo Transmission Systems}

To assess the impact of gas composition on network performance, we simulated the same transport scenario with varying levels of CO$_2$ purity using MYNTS. The flow conditions and boundary parameters were kept identical to those used in the base case. Three different mixtures were considered; a pure case with 100\% CO$_2$, a slightly impure case with 98\% CO$_2$, 1\% H$_2$, 1\% O$_2$, and a more impure case with 96\% CO$_2$, 2\% H$_2$, and 2\% O$_2$. These impurity combinations were selected to represent plausible compositions of captured CO$_2$ streams that could arise from industrial or power plant processes. As discussed by Aursand et al. \cite{aursand2013pipeline}, CO$_2$ transported in carbon capture and storage systems is often not entirely pure and may include small amounts of impurities such as O$_2$, H$_2$, and N$_2$, depending on the capture technology and fuel source.

Initially, the same pump configurations as in the 100\% CO$_2$ case were applied to the other cases. However, due to the altered thermophysical properties of the impure mixtures, these settings did not maintain the desired operating conditions. Therefore, the power and outlet temperature of pumps were iteratively adjusted to achieve stable operation in each case.

\begin{figure}[H]
    \centering
    \includegraphics[width=\linewidth]{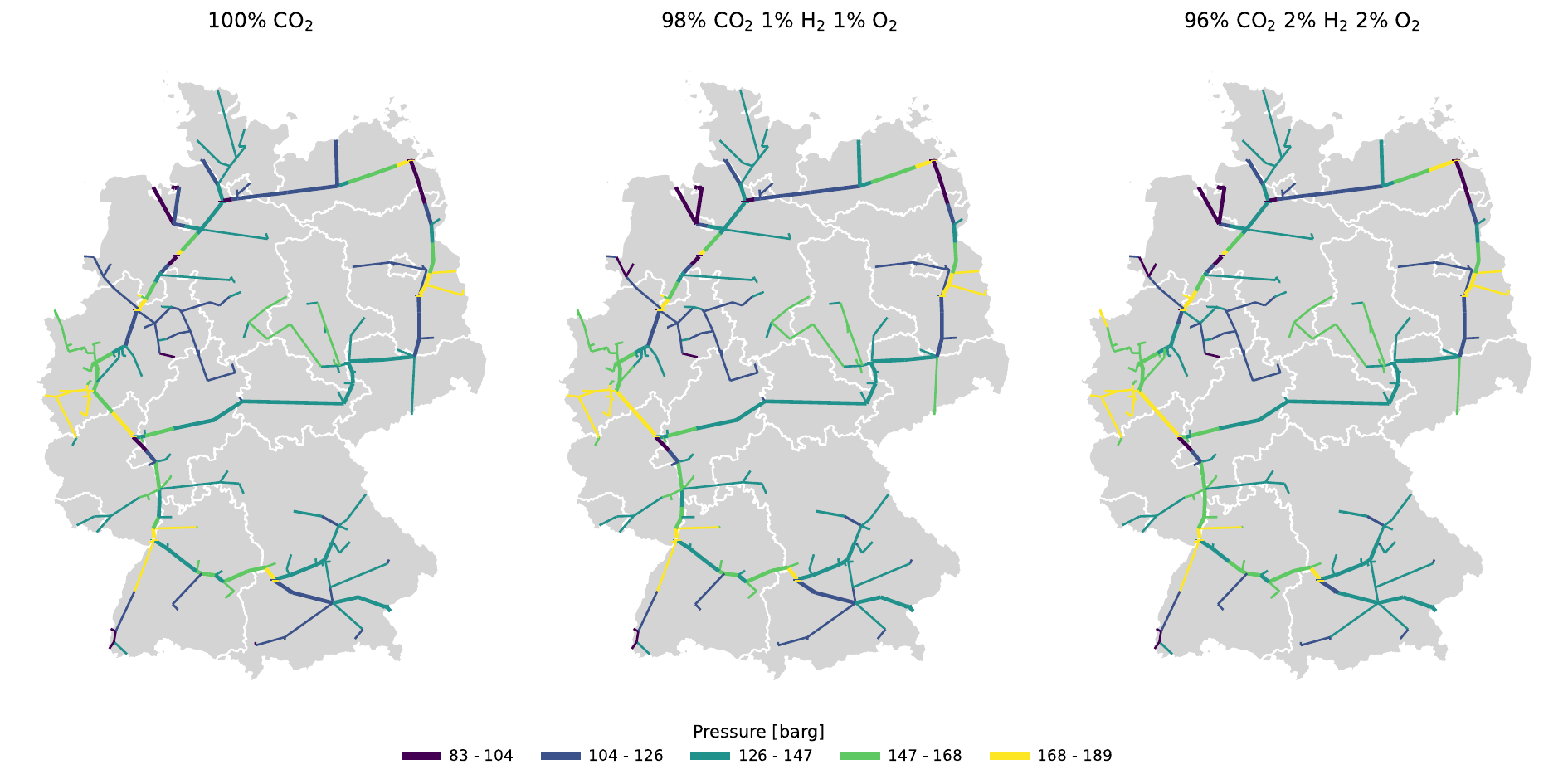}
    \caption{Pressure results of CO$_2$ network simulations with different impurity levels. The three cases correspond to 100\% CO$_2$ (left), 98\% CO$_2$, 1\% H$_2$, 1\% O$_2$ (middle), and 96\% CO$_2$, 2\% H$_2$, 2\% O$_2$ (right).}
    \label{fig:impurities}
\end{figure}

Figure~\ref{fig:impurities} illustrates the resulting pressure distributions. Despite compositional differences, all three cases remain within a reasonable pressure–temperature range, ensuring stable dense phase transport for the tested impurity levels. This demonstrates that the network can accommodate compositional variations in CO$_2$ streams with only moderate operational adjustments.

It should be noted that the simulation of CO$_2$ mixtures containing impurities in most cases will slow the simulation speed. Therefore, it is recommended to start the analysis of pump configurations and other optimization steps always with pure CO$_2$. One could add a kind of "safety distance" of 5 bar or more above the critical pressure to avoid unwanted phase transitions. However, for the final configuration, a simulation with impurities is recommended.

\section{Discussion}

This study has developed a method for comprehensive assessment of future CO$_2$ pipeline systems by integrating energy systems modeling with physical pipeline simulation. The combination out of energy system scenario data (Langfristzenarien III) and physical simulation with MINTH bridges a research gap between the system-focused approach of energy system modeling and the physical, engineering-based modeling of pipeline systems. 
The strengths of this methodology are:
\begin{itemize}
    \item The use of projected future CO$_2$ emissions based on energy system scenarios, which are aligned with climate targets and economic development.
    \item The site-specific resolution, which specifies actual distances, thereby enables concrete cost estimates and also relates to actual industrial facilities and their emissions.
    \item Aligning pipeline routes with existing infrastructure (the natural gas network), which offers advantages during implementation.
    \item Detailed analysis of the physical conditions within the pipelines, taking into account impurities, elevation differences, temperature and pipeline diameters.
    \item Determining locations for pumps to ensure safe operation in dense phase. 
    \item Estimating system costs.
\end{itemize}
the use of projected future CO2 emissions based on energy system scenarios, which are thus aligned with climate targets.
The site-specific resolution, which specifies actual distances, thereby enables concrete cost estimates and also relates to actual industrial facilities and their emissions.
Aligning pipeline routes with existing infrastructure (the gas network), which offers advantages during implementation
Detailed analysis of the physical conditions within the pipelines, taking into account impurities, elevation differences, and temperature.
Determining locations for pumps to ensure safe operation during the compression phase
The method was applied to a case study for Germany. It examined the design and characteristics of a prospective CO$_2$ pipeline network for the year 2045. While the analysis provides insights into the spatial structure and potential scale of such a network and validates the method developed, the single topology and scenario calculated is not a sufficient basis for infrastructure planning. In the following, we discuss ways forward to extend the analysis towards a more comprehensive basis for infrastructure planning and carbon management strategy development.

In this case study, we assumed that offshore storage capacity has no limitations, neither with regard to total sequestration volume nor with regard to sequestration rates. Storage availability can, however, become an important bottleneck. At present, there are no operational storage projects in Germany, and it remains uncertain whether the capacities of already announced storage projects in the EU will be sufficient to accommodate the projected volumes \cite{InternationalAssociationofOil&GasProducers.2025}. The geological storage potential appears to be adequate, and this also applies to large parts of the German onshore area, but the actual deployment is uncertain. In light of the recent Carbon Storage and Transport Act (KSpTG), which enables onshore storage, future studies should explicitly include onshore storage options in Germany and assess the storage capacities at a European scale \cite{FederalGovermentGermany.2025}.

While we have assessed one set-up CO$_2$ sources and sinks, the high uncertainty regarding the future locations and quantities should be reflected by additional scenarios. E.g., the network configuration for 2045 is assumed to follow the spatial distribution of industrial sites as observed today. Plant closures or changes in capacity can substantially alter this assumption and relocation of green industrial processes might change entire value chains, e.g., in the chemical industry \cite{Verpoort.2024}. Further, we have assessed a network for limited carbon sources from cement and lime production and waste incineration. More industries could qualify for carbon capture with high volume point sources like steel production, steam crackers, steam reformers and large chemical industry clusters or even power plants. Over the longer term, carbon capture to enable negative emissions also needs to be considered. If high levels of demand for negative emissions materialize, it can be expected that additional biomass-based facilities for generating biogenic negative emissions will be established, and that existing plants will be connected to the network. Consequently, both the number and location of emission sources and sinks, as well as their emission or uptake volumes, are likely to change over time. The pipeline network will therefore need to be designed with sufficient flexibility to accommodate such structural changes. This flexibility is not yet reflected in the topology analysed here and should be subject to future analyses.

To determine whether CCUS is an economically viable alternative to other emission abatement strategies, a holistic system perspective is required. The transport volumes used in this study are derived from Energy system and infrastructure models, which are well suited for this purpose. Nevertheless, there is scope for methodological improvement by integrating pipeline infrastructure directly into the optimization framework and allowing the model to decide, for each site, whether and how it should be connected to the network. Hofmann et al.~\cite{Hofmann.2025} have demonstrated such an approach in a simplified form. Extending this type of analysis to a spatial resolution comparable to that used here, and coupling it with detailed pipeline design, could provide more accurate and more robust results.

In this study, the network planning was carried out under the assumption of transport in the dense phase. Under these conditions, it is unlikely that existing pipelines can be reused directly, as they were not designed for this aggregate state and associated pressure regimes. In contrast, the gaseous phase can be used for transporting smaller quantities of CO$_2$. This mode imposes less stringent requirements on the pipeline itself but requires more frequent compression over longer distances. As a result, gaseous transport may be particularly suitable for short collection networks feeding into centralized liquefaction hubs, from which the CO$_2$ is then transported in a main dense-phase pipeline. Moreover, gaseous transport can offer advantages over high-pressure dense-phase pipelines in situations where stricter safety requirements apply. Identifying where conversions of existing pipelines to gaseous CO$_2$ transport are technically and economically feasible requires detailed information about the condition and properties of the existing assets, data that are usually not publicly available. Nevertheless, it can be assumed that such measures could lead to cost savings in specific cases.

Another simplification in the present analysis is the assumption that the CO$_2$ pipeline routes follow the existing natural gas network and that industrial sites are connected via direct lines. Routing CO$_2$ pipelines along existing gas corridors offers clear advantages, especially in terms of reduced planning and permitting times. However, other factors can significantly influence the optimal route. For dense-phase transport, CO$_2$ behaves as a liquid, so a route with minimal elevation changes is desirable in order to avoid unnecessary pumping and associated energy consumption. Such influences can be incorporated using a least-cost path approach that assigns cost factors to different geographic and technical constraints. Yeates et al.~\cite{Yeates.2024} provide an example of such a methodology. Applying a similar approach would allow the connection lines, which are assumed to be direct in this study, to be represented more realistically. Based on their least-cost routing, Yeates et al.\ obtain a total network length of between 4300 and 5200~km, compared with approximately 7000~km in the topology developed here. However, this comparison should be treated with caution, as the underlying sets of sources and sinks differ; for example, Yeates et al.\ also include emissions from the power sector. Nonetheless, their approach can be considered more accurate with respect to spatial routing and should be taken into account in future analyses.

Finally, this study has focused on pipeline transport, but alternative transport modes such as ship or rail also represent viable options in a comprehensive CO$_2$ infrastructure system. Due to their different cost structure, characterized by higher operational expenditures but lower capital costs compared to pipelines, these modes are particularly attractive for small to medium transport volumes and for flexible or temporary connections. They should therefore be included in future assessments of integrated CCUS infrastructure, especially for transitional time horizon. To date, comparisons of different transport options have mainly been carried out in case studies, such as \cite{Oeuvray.2024} and \cite{Becattini.2024}. Extending such comparative analyses to system-scale studies for Germany and Europe would provide a more complete basis for evaluating the role of CCUS in long-term decarbonization strategies.


\subsection{Related work}
In this study, we combined a network optimization tool with a physical simulator. This approach is - to some extent - similar to the techniques used in SimCCS \cite{Middleton.2009, Gunawan.2024} with some important differences: In SimCCS, the model for CO$_2$ transport has been strongly simplified to integrate it into the optimization task. Here, we use a physical simulation tool to predict the behavior of the transported fluid in a suitable way for network planning including topographic information, pressures, phase transitions, etc. The drawback of our approach is the fact that this simulator is not integrated into the route optimizer. This means that an iterative loop has to be performed which might be more time consuming compared to SimCCS but will be more accurate for the final results in a physical sense.

\section{Conclusion and Outlook}

This paper has presented an integrated approach to planning a national CO$_2$ transport infrastructure by combining long-term energy system scenarios with detailed physical network simulation. We calculated spatially explicit CO$_2$ balances for a scenario that includes carbon capture from hard-to-abate sources from cement and lime production as well as waste incineration. It considers substantial cross-border transit imports of CO$_2$ from neighboring countries and four coastal hubs for offshore-storage in the North Sea, but no onshore storage nor carbon use.

On this basis, we developed a CO$_2$ pipeline topology that connects all relevant components. The topology was constructed with high spatial resolution and follows existing natural gas transmission corridors wherever possible, in order to reduce planning and permitting efforts. The resulting backbone network has a total length of approximately 7000~km and is designed as a new-build system for dense-phase CO$_2$ transport.

Using the multiphysical network simulator MYNTS, we then assessed the technical feasibility and optimized the physical design of this network. We evaluated uniform and hybrid diameter configurations and found that a flow-sensitive design combining DN700 pipelines in high-flow backbone segments with DN500 and DN400 pipelines in low-flow branches can achieve similar total pipeline investment costs to a uniform DN500 system, while significantly reducing the number of pumps required. This configuration is associated with total pipeline investments on the order of 17~billion~€.

The simulations further showed that elevation differences have a non-trivial impact on pressure distribution in dense-phase CO$_2$ networks. The hydrostatic pressure contribution can lead to situations where downstream nodes at lower elevation exhibit higher pressures than upstream nodes, so neglecting terrain height would underestimate pressure variations and could lead to misplacement or misdimensioning of equipment. In addition, the analysis demonstrated that, for the impurity levels considered here (up to 4\% admixture of H$_2$ and O$_2$), the network can maintain stable dense-phase operation with moderate adjustments to pump settings and cooling, provided that pressure and temperature stay within the dense-phase envelope. Phase diagrams for all network nodes confirmed that operating points remained safely above the gas–dense-phase boundary.

Overall, the combination of scenario-based system analysis, spatial topology development, and detailed physical flow simulation proofs valuable for CO$_2$ infrastructure planning. It captures operational and thermodynamic constraints—such as pressure maintenance, elevation effects, velocity limits, and phase behavior—that are not represented in high-level system models but are essential for realistic cost and feasibility assessments. This integrated approach allows more reliable statements about the design and dimensioning of future CO$_2$ transport systems and can inform both strategic planning and early-stage project development. The method developed is similarly applicable to other countries and regions.\\

While the defined scenario and case study for Germany validates the method, the single topology and scenario calculated is not a sufficient basis for infrastructure planning. The case study should be extended in three main directions. First, a European perspective should be added, explicitly modeling cross-border pipeline links, ship-based connections, and offshore storage hubs. This would allow a consistent assessment of transit flows through Germany and the role of shared storage capacities in countries around the North Sea. Second, the static 2045 view should be complemented by a time-resolved analysis of the ramp-up phase. This includes the gradual deployment of the backbone, evolving spatial patterns of CO$_2$ sources and sinks, and the temporary role of smaller-scale or vehicle-based transport (ship and rail) before a fully developed pipeline network is in place. Third, more uncertainties should be tested with additional scenarios and topologies. This includes the potential role of onshore storage in Germany, the possible reuse or partial conversion of existing pipelines and corridors, and higher shares of negative emissions and CO$_2$ utilization as well as other CO$_2$ capture in other industries like steel or chemicals. 


\section{Acknowledgments}
We gratefully acknowledge the partly support of the German Federal Ministry of Research, Technology, and Space in the project
TransHyDE-Sys, grant 03HY201M. as well as from the Fraunhofer Cluster
CINES, and from the EU HORIZON project \emph{Eastern Lights}.

\section{Declarations}
\subsubsection*{Ethics approval and consent to participate}
Not applictable
\subsubsection*{Consent of publication}
The authors have agreed to the publication. Other individuals are not affected.
\subsubsection*{Competing interests}
The authors declare no competing interests. 
\subsubsection*{Resarch data for this article}
See appendix. Data or code beyond this cannot be provided.
\subsubsection*{Authors contribution}
\textbf{M.A.}: Conceptualization, Methodology, Software, Investigation,  Visualization, Supervision, Project administration, Funding acquisition
\textbf{M.N.}: Conceptualization, Methodology, Software, Investigation, Writing – Original Draft, Writing – Review \& Editing, Visualization, 
\textbf{O.A.}: Conceptualization, Software, Investigation, Writing – Original Draft, Writing – Review \& Editing, Visualization
\textbf{L.L.}: Methodology, Investigation, Writing – Original Draft, Writing – Review \& Editing, Visualization, 
\textbf{S.B.}:  Methodology, Investigation, Writing – Original Draft, Writing – Review \& Editing, Visualization,
\textbf{T.F.}: Conceptualization, Methodology, Writing – Review \& Editing, Visualization, Supervision, Project administration, Funding acquisition
\textbf{B.K.}: Conceptualization, Methodology, Investigation, Writing – Review \& Editing, Visualization, Supervision, Funding acquisition


\section*{Appendix}

The CO$_2$ quantities used per location can be found in \href{appendix/location_data_HtA+Imports}{location data HtA+Imports.csv}.
Node and element lists are provided to document the network's topology.

\cleardoublepage
\bibliographystyle{plain}   


\end{document}